\documentclass[12pt]{article}
\usepackage {amsmath}
\input amssym.def
\input amssym
\input{psfig} 
\def \C{{\Bbb C}}

\def \R{{\Bbb R}}

\def \L{{\cal L}}
\def \l{{\lambda}}

\def \Z{{\Bbb Z}}

\def \BC{{\Bbb C}}

\def \kl {K_G^*(X)}

\def \omg {\omega}
\def \oo {\Omega}
\def \tl {\sim}

\def \un{\underline}
\def \BZ{\Bbb Z}
\def \BQ{\Bbb Q}
\def \BC{\Bbb C}
\def \hook{\hookrightarrow}
\def \tt {{\cal T}}
\def \proof{{\noindent{\it Proof.\ \ }}}
\newtheorem{Th}{THEOREM}[section]

\newtheorem{prop}[Th]{PROPOSITION}

\newtheorem{lemma}[Th]{LEMMA}

\begin{document}
\title{Equivariant K-Theory of Simply Connected Lie Groups}
\author {Jean-Luc Brylinski\\
 Department of Mathematics\\The  Pennsylvania State University \\
University Park, PA 16802, USA\\
        {\sf email: jlb@math.psu.edu } \\
Bin Zhang\\
 Department of Mathematics\\The  Pennsylvania State University \\
University Park, PA 16802, USA\\
        {\sf email: zhang\_b@math.psu.edu } }

\maketitle

\begin{abstract}
We compute the equivariant $K$-theory $K_G^*(G)$ for a simply connected Lie group $G$ (acting on itself by conjugation). We prove that $K_G^*(G)$ is isomorphic to the algebra of Grothendieck differentials  on the representation ring.  We also study a special example of a non-simply connected Lie group $G$, namely $PSU(3)$, and compute the corresponding equivariant $K$-theory.
\end{abstract}

\section {Introduction}

For a finite group $G$, the equivariant $K$-theory of the $G$-space $G$ (where the $G$-action is the conjugation action) has a very nice expression \cite{gl1}:
$$K_G^0(G)\otimes \C=\{{\rm \ conjugacy\  invariant\  functions}\  {\cal C}\to \C\}$$
where ${\cal C}$ is the set of $\{(g,h)\in G \times G | gh=hg\}$. 

Since the group $SL(2,\Z)$ acts on the set ${\cal C}$, we have an action
of $SL(2,\Z)$ on the group $K^0_G(G)$. In particular, the element
$\left (\begin {array}{ll}0&1\\ -1&0\end{array}\right)$ gives the so-called Fourier transformation,
 which was used by Lusztig \cite {gl1} in the representation theory  of finite reductive groups. The space $K^0_G(G)$ appears as the space of conformal
field theories attached to a finite group \cite {dvvv}.  So it is natural
to ask what happens  when $G$ is a compact group.

The first author conjectured that for a compact group $G$, the equivariant
K-theory $K^*_G(G)$ is isomorphic to the algebra
$\oo ^*_{R(G)/\Z}$ of Grothendieck differentials, where $R(G)$ is the 
representation ring of $G$. 

This was motivated by his computation
of the equivariant cohomology theory $H^*_G(G,\C)$ of $G$-space $G$ with adjoint action.

For a compact group $G$ and a compact $G$-space $X$, the equivariant cohomology $H^*_G(X)$ is defined as $H^*_G(X)=H^*(EG{\overset G \times}X)$ \cite{ma1}.

If we view $G$ as a $G$-space where the action is conjugation, we can prove that $H_G^*(G)\cong \Omega ^*_{H^*(BG,\C)/\C}$ (with coefficients $\C$). 
There are several proofs of this result, one of which is based on the fact
that $H_G^*(G)$ is isomorphic to $H^*(B(\L G))$, where $\L G$ is the smooth
loop group of $G$.

The main result of the paper is the construction of an isomorphism
$\oo ^*_{R(G)/\Z}\simeq K^*_G(G)$ in the case where $G$ is simply connected
(Theorem 3.2). Though we always have an algebra map
$\phi:\oo ^*_{R(G)/\Z}\to K^*_G(G)$ for any compact group $G$, we show that in case
 $G=PSU(3)$, $\phi$ is not an isomorphism. In fact 
in that case $K^*_G(G)$ has torsion as a $R(G)$-module

For a simply-connected Lie group $G$, the proof of our result is as follows. By using  Hodgkin's spectral squence, we describe the $R(G)$-module structure of the equivariant $K$-theory $K_G^*(G)$, we prove that is a free module over the representation ring $R(G)$. Then by considering the natural map from $\Omega _{R(G)/\Z}^*$ to $K_G^*(G)$,  we obtain the algebra structure.

{\bf Thanks} This research was supported in part by NSF grant DMS-9504522 and DMS-9803593.

\section{Background}

In this section we recall some basic material we will use later: representation rings, equivariant $K$-theory and the algebra of Grothendieck differentials. Here we always assume $G$ is a compact Lie group.

\subsection {The representation ring}

For a compact Lie group $G$, a $G$-module means a finite-dimensional complex vector space $M$ with a continous linear $G$-action on it. If $M$ and $N$ are $G$-modules, we can form their direct sum $M\oplus N$, and with respect to this operation, the isomorphism classes of $G$-modules form an abelian semigroup; the associated abelian group $R(G)$ is called the representation ring of $G$, the tensor product induces a commutative ring structure in $R(G)$.

A $G$-module $M$ has a character: $\chi _M:\ G \rightarrow \BC, \ \chi _M(g)=tr _M (g)$, and $M$ is determined by $\chi _M$ up to isomorphism. We know
$ \chi _{M \oplus N}=\chi _M+\chi _N$, $ \chi _{M \otimes N}=\chi _M\chi _N$ and $\chi _M(hgh^{-1})=\chi _M (g)$, thus the map $M\mapsto \chi _M$ identifies $R(G)$ with a subring of the ring of $G$-invariant complex functions on $G$, so $R(G)$ is also called the character ring of $G$. 

We already know \cite{gs2} that the representation ring of the torus $T^n$ is the ring of Laurent polynomials,
$$R(T^n)={\Bbb Z}[X_1,X_2, \cdots, X_n, (X_1X_2 \cdots X_n)^{-1}]$$

The following result is critical to us:

\begin {lemma} {\rm \cite {nb1}} For a simply connected Lie group $G$, $R(G)$ is 
a polynomial ring over $\Z$ with $rank(G)$ generators.
\end {lemma}

\subsection {Equivariant $K$-theory}

Let $X$ be a locally compact $G$-space, then the equivariant $K$-theory $K_G^*(X)$ is defined and it is a $\Z _2$-graded algebra \cite{gs1}. In case that $X$ is compact, $K_G^0(X)$ is generated by equivalence classes of complex $G$-vector bundles over $X$.

Equivariant $K$-theory is a common generalization of  two
important extreme cases.
When  $G=1$, $K_G^*(X)$ is reduced to the ring $K^*(X)$ of the
ordinary $K$-theory. 
When
$X$ is a single point, $K_G^0(X)$ is nothing  but the representation
ring $R(G)$ and $K^1_G(X)=0$.  

$K_G^*(X)$ enjoys some very nice properties:

1. $K_G^*$ is functorial with respect to both the proper 
equivariant maps between spaces  
and the homomorphisms between groups. So, for any compact $G$-space $X$, 
the obvious map $X \to {\rm point}$ gives rise to a $R(G)$-module structure on $K_G^*(X)$, therefore $R(G)$ serves as the coefficient ring 
in equivariant $K$-theory.

If $B$ is a closed subgroup of $G$, for any $B$-space $X$, we get the map $i^* : K_G^*(G {\overset B \times} X) \to K_B^*(X)$ induced from the map $i: B \hook G$, which in fact is an isomorphism.

As a special case, we have:

$$K_G^q(G/B)=K_B^q({\rm point})=\left\{ \begin{array}{cc}
R(B), & q = 0 \\
0, &  q=1 
\end{array} \right.
$$

2. Bott periodicity:

The Thom homomorphism \cite {gs1} $\phi_*: K_G^*(X) \rightarrow K_G^*(E)$ is an isomorphism for any complex $G$-vector bundle $E$ over $X$.

3. Six term exact sequence: for any pair $(X,A)$, where $A$ is closed $G$-subset in $X$, we have the exact sequence:

$$\begin{array}{ccccc}
K_G^{-1}(X,A)&\rightarrow &K_G^{-1}(X) &\rightarrow & K_G^{-1}(A)  \\
\uparrow&&&&\downarrow \\
K_G^{0}(A)&\leftarrow &K_G^{0}(X) &\leftarrow &K_G^0(X,A)
\end{array} $$

Also we have an effective way to compute the equivariant $K$-theory in many cases, that is the Segal spectral sequence:

For any covering $\cal U$ of $X$ by closed $G$-stable subsets, there exists a multiplicative spectral sequence $E^{*,*}_n \Rightarrow \kl$ \cite{gs3}, such that $E_2^{p,q}={\check H} ^p({\cal U}, {\cal K}_G^q)$, where ${\cal K}_G^q$ means the coefficient system: $U \mapsto K_G^q(U)$.

\subsection {Grothendieck Differentials}

Let $A \subset B$ be commutative rings, the algebra of Grothendieck differentials $\Omega ^*_{B/A}$ \cite{ag1} is the differential graded $A$-algebra constructed as follows:

Let $F$ be the free $B$-module generated by all elements in B, to be clear, we use $db$ to denote the generator corresponding to $b \in B$, so 
$$F=\bigoplus _{b\in B} Bdb.$$
and let $I \subset F$ be the $B$-submodule generated by 
$$\left\{ \begin{array}{l}
da,\ \forall a \in A\\
d(b_1+b_2)-db_1-db_2, \ \forall b_1, b_2 \in B\\ 
d(b_1b_2)-b_1db_2-b_2db_1,\ \forall b_1, b_2 \in B \end{array} \right\},$$ we then get the quotient $B$-module $$\Omega _{B/A}=F/I.$$

Let $\Omega ^0_{B/A}=B$, $\Omega ^1_{B/A}=\Omega _{B/A}$, and $\Omega ^p_{B/A}= \Lambda ^p_B \Omega _{B/A}$, there is a differential: $d: \Omega ^p_{B/A} \to \Omega ^{p+1}_{B/A}$, which maps $b \in B$ to $db$, then  
$$\Omega ^*_{B/A}=\bigoplus _{p=0}^{\infty} \Omega ^p_{B/A} $$ 
is the differential graded algebra of Grothendieck differentials of $B$ over $A$. It is the generalization of the algebra of differentials on affine spaces, for example, if $B=A[x_1, \cdots, x_n], $, then $\Omega ^p_{A[x_1, \cdots, x_n]/A}=\oplus _{i_1<i_2< \cdots <i_p} A[x_1, \cdots, x_n]dx_{i_1}\wedge \cdots \wedge dx_{i_p}$.

\section {Construction of $\phi$ and main result}

In this section, we construct the map $\phi$ from $\oo ^*_{R(G)/\Z}$ to $K^*_G(G)$
and state our main result.

We can define the map $\phi$ from $\Omega ^*_{R/\Z}$ to $K_G^*(G)$ ($G$ acts on itself by conjugation) as follows:

For any representation $\rho : G\to GL(V)$, we let $\phi (\rho)$ be the 
 $G$-vector bundle $G\times V$ over $G$.
 Then let $\phi (d\rho) \in K^1_G(G)$ be given by the complex of $G$-bundles over $G \times \R$:

$$\begin {array}{ccccc}
0\to &G\times \R \times V&\to &G\times \R \times V&\to 0\\
&(g, t, v)&\mapsto&(g,t,t\rho (g)v)&
\end{array}
$$
This complex is exact out of $G\times \{0\}$, 
thus it defines an element in $K_G^{1}(G)$. 
By the properties of equivariant $K$-theory, we know $\phi (d\rho ) \cup \phi ( d\rho)=0$.
 
\begin {prop} There is a unique algebra homomorphism $\phi: \Omega ^*_{R(G)/\Z} \to K^*_G(G)$ for which $\phi (\rho)$ and $\phi (d\rho)$ are as described above.
\end{prop}

\proof By the definition of $\phi$, we only need to prove that for any $\rho _i: G \to GL(V_i)$, $i=1,2$, $\phi (d(\rho_1 \rho_2))= \phi (\rho_1)\phi ( d\rho_2)+\phi (\rho_2)\phi ( d\rho_1)$ in $K_G^{1}(G)$. Once we prove this, then we can extend $\phi$ to an algebra  map from $\Omega ^*_{R/\Z}$ to $K_G^*(G)$.   

$\phi (d(\rho_1 \rho_2))$ is the complex:
$$\begin {array}{ccccc}
0\to &G\times \R \times (V_1\otimes V_2)&\to &G\times \R \times (V_1\otimes V_2)&\to 0\\
&(g, t, v_1\otimes v_2)&\mapsto&(g,t,t\rho_1(g)(v_1) \otimes \rho_2(g)(v_2))&
\end{array},
$$

$\phi (\rho_1)\phi( d\rho_2)$ is the complex:
$$\begin {array}{ccccc}
0\to &G\times \R \times (V_1\otimes V_2)&\to &G\times \R \times (V_1\otimes V_2)&\to 0\\
&(g, t, v_1\otimes v_2)&\mapsto&(g,t,tv_1\otimes \rho_2 (g)(v_2))&
\end{array}
$$

And
$\phi (\rho_2)\phi ( d\rho_1)$ is the complex:
$$\begin {array}{ccccc}
0\to &G\times \R \times (V_1\otimes V_2)&\to &G\times \R \times (V_1\otimes V_2)&\to 0\\
&(g, t, v_1\otimes v_2)&\mapsto&(g,t,t\rho_1(g)(v_1)\otimes v_2)&
\end{array}
$$
Now we invoke the following fact:

The two maps $GL(n)\times GL(n)\to GL(2n)$ given by
$$(A,B)\mapsto \left(\begin{array}{cc}
               A&0\\
               0&B
\end{array}\right)
$$
and
$$(A,B)\mapsto \left(\begin{array}{cc}
               AB&0\\
               0&1
\end{array}\right)
$$
are homotopic. The homotopy is given explicitly by 
$$\rho _s (A,B)=
\left (\begin {array}{cc}
A&0\\
0&1
\end{array}\right )
\left (\begin {array}{cc}
\cos s&\sin s\\
-\sin s&\cos s
\end{array}\right )
\left (\begin {array}{cc}
1&0\\
0&B
\end{array}\right )
\left (\begin {array}{cc}
\cos s&-\sin s\\
\sin s&\cos s
\end{array}\right )
$$
where $0\le s \le \pi /2$.
It is $GL(n)$-equivariant in the following sense:
$$
\rho _s (g A g^{-1}, g B g^{-1})= \left (\begin {array}{cc}
g&0\\
0&g
\end{array}\right )\rho _s (A, B)
\left (\begin {array}{cc}
g&0\\
0&g
\end{array}\right )
^{-1}$$

Following the same construction of homotopy, we can prove
$$\phi (d(\rho_1 \rho_2))\oplus E\simeq \phi (\rho_1)\phi ( d\rho_2)+\phi (\rho_2) \phi ( d\rho_1)$$
where $E$ is the trivial complex:
$$0 \to G\times \R \times (V_1\otimes V_2) \to G\times \R \times (V_1\otimes V_2) \to 0 $$
thus
$$\phi (d(\rho_1 \rho_2))= \phi (\rho_1) \phi (d\rho_2)+\phi (\rho_2) \phi ( d\rho_1)$$ in $K_G^{1}(G)$. Q.E.D.

The following theorem is the main result of this paper and will be proved in sections 4-6.

\begin {Th} For a compact, simply connected Lie group $G$, the algebra homophism $\phi:\ \Omega ^*_{R/\Z} \to K^*_G(G)$ is an algebra isomorphism.
\end {Th}

\section {Module structure}

From now on we assume $G$ to be a simply connected Lie group. $T$ is a maximal torus of $G$, $W$ is the corresponding Weyl group. In this case, we know by Lemma 2.1 that the representation ring $R$ of $G$ is a polynomial ring. As usual, we make $G$ act on itself by conjugation. 

\begin {lemma} $K_G^*(G)\cong K^*_{G\times G}(G\times G)$
where the action of $G\times G$ on $G\times G$ is given by: $(g_1, g_2).(h_1,h_2)=(g_1h_1g_2^{-1}, g_1h_2g_2^{-1})$
\end {lemma}

\proof Consider $G$ as the diagonal subgroup of $G\times G$, the functoriality of equivariant $K$-theory yields, 
$$K_G^*(G)\cong K^*_{G\times G}({(G\times G)\overset {G\times G}\times}G)$$
where we view $G$ as a $G\times G$-space by the following action:
$(g_1, g_2).(h_1)=g_1h_1g_2^{-1}$. Then we have a $G \times G$-equivariant homeomorphism:
$$(G\times G){\overset {G}\times}G\cong G\times G.$$ 
where the homeomorphism is:
$$ ( (g_2, g_3),g_1)\mapsto (g_2 g_1 g_3^{-1}, g_2 g_3 ^{-1})$$
Obviously this map is a $G\times G$-map. 
Q.E.D.

Let us recall a theorem of Hodgkin here:
 
\begin {Th} {\rm \cite {lh} \cite {gl2}} For a Lie group $G$, supposing $\pi _1(G)$ is torsion free, then  for any $G$-spaces $X$, $Y$ there is a spectral sequence $E_r\Rightarrow K_G^*(X\times Y)$, with $E_2$-term

$$E_2^{*,*}=Tor_{R(G)}^{*,*}(K_G^*(X), K_G^*(Y))$$
\end{Th}

\underbar{\it Remark on $Tor ^{*,*}_R$.} Let $R$ be a commutative ring, $M^{\bullet}$, $N^{\bullet}$ be graded differential complexes over $R$, then the bi-complex $Tor_{R}^{*,*}(M^{\bullet}, N^{\bullet})$ is defined (for example, see \cite {jm1}).
If we take a proper projective resolution $P^{\bullet}=\{ (P^r)^s \}$ of $N^{\bullet}$, then $Tor_{R}^{*,*}(M^{\bullet}, N^{\bullet})$ is the homology of the total complex of $M^{\bullet}\otimes Total (P^{\bullet})$. It is bigraded:
$$\begin {array}{l}
H^i(M^{\bullet}\otimes _R Total (P^{\bullet}))\\
={\rm subquotient\ of\ } \oplus _{m+n=i} M^m\otimes _R Total (P^{\bullet})^n\\
={\rm subquotient\ of\ }\oplus _{m+n=i} M^m\otimes _R (\oplus_{r+s=n} (P^r)^s)
\end{array}$$
Those elements in it having degree $m+s=j$ form $Tor^{i,j}_R(M^{\bullet},N^{\bullet})$.
 
In our case, if we view $G$ as a $G\times G$-space (under the above action), then 
$$K^i_{G\times G}(G)=\left \{ \begin{array}{ll}
R,&i=0\\
0,&i=1
\end{array}\right .$$
and we find that there is a spectral sequence $E_r\Rightarrow K_G^*(G)$ such that 
$$E_2^{*,*}=Tor_{R\otimes R}^{*,*}(R,R)$$
If consider the bidegree, we see that $M^{0}$=$N^{0}$=$R$, and $M^{i}$=$N^{i}=0$ for $i>0$, so $Tor^{i,j}_R(M^{\bullet},N^{\bullet})$ can be non-zero only when $j=0$.

\begin {lemma} We have $Tor_{R\otimes R}^{*,*}(R,R)\cong \Omega ^*_{R/\Z}$. In particular, $Tor_{R\otimes R}^{*,*}(R,R)$ is a free $R$-module of rank $2^{rank(G)}$.
\end{lemma}

\proof Because $R$ is a polynomial ring, $R=\Z[x_1, \cdots, x_n]$, we have
$$R\otimes R=\Z[x_1, \cdots, x_n,y_1, \cdots, y_n].$$ 
By using the Koszul 
complex we can get a free resolution of $R$ over $R\otimes R$:
$$
\cdots \to\oplus_{i_1 <\cdots < i_r} \Z[x_1, \cdots, x_n,y_1, \cdots, y_n]e_{i_1}\wedge \cdots \wedge e_{i_r}{\overset {e_i \mapsto x_i-y_i} \longrightarrow} 
$$
$$\oplus _{i_1 <\cdots < i_{r-1}} \Z[x_1, \cdots, x_n,y_1, \cdots, y_n]e_{i_1}\wedge \cdots \wedge e_{i_{r-1}}\to
$$
$$\cdots \to \Z[x_1, \cdots, x_n,y_1, \cdots, y_n] \to \Z[x_1, \cdots, x_n]
$$
  Using this resolution we can prove the statement easily. Q.E.D.

\begin{lemma} $rank_R(K_G^*(G))=2^{rank(G)}$
\end{lemma}

\proof Firstly by the localization theorem \cite{gs1}
 we know that
 $$rank_{R(T)}K^*_T(G)=rank_{R(T)}K^*_T(G^T)=rank_{R(T)}K^*_T(T)=2^{rank(G)},$$
where $G^T$ denotes the fixed point set of $T$ in $G$.
 
Secondly, by \cite {gl2} \cite {lh}, we know that
 $$K_T^*(G)=K^*_G(G)\otimes _{R(G)} R(T)$$
 So the rank of $K_G^*(G)$ over $R$ is the same as the rank of $K_T^*(G)$ over $R(T)$. Q.E.D.

\begin{lemma} We have $K_G^*(G)\cong \Omega _{R/\Z}^*$ as $R$-module.
\end{lemma}

\proof $E_2^{*,*}=Tor_{R\otimes R}^{*,*}(R,R)\cong \Omega ^*_{R/\Z}$ is a free $R$-module of rank $2^{rank(G)}$, the spectral sequence converges to $K_G^*(G)$, which has the same rank as $E_2$, hence $d_r, \ r\ge 2\ =0$. Since $E_2$ is a free $R$-module, $E_{\infty}$ is isomorphic to $E_2$ as an $R$-module.
Q.E.D.

\section {The Algebra Structure}

So far we have discussed the $R(G)$-module structure on $K^*_G(G)$. Let us now consider the 
algebra structure. 

\begin {lemma} $K^*_T(T)\cong \Omega ^*_{R(T)/\Z}$
\end{lemma}

We will prove the following lemma in the next section.

\begin {lemma} The inclusion $j: R(G)\to R(T)$ induces an isomorphism: 
$j^* : \Omega ^*_{R(G)/\Z} \to (\Omega ^*_{R(T)/\Z})^W$.
\end{lemma}

We already constructed the algebra homomorphism  $\phi :\Omega ^*_{R/\Z} \to 
K^*_G(G)$ in section 3. Now we can prove

\begin{Th} $\phi$ is an isomorphism.
\end {Th}
\proof Let us consider the maps:
$$\Omega ^*_{R/\Z}{\overset {\phi} \to} K^*_G(G){\overset {\alpha} \to} K^*_T(G){\overset {\beta} \to}K^*_T(T) \cong \Omega ^*_{R(T)/\Z}$$
where $\alpha$ and $\beta$ are all natural maps.

The composition $\beta\circ \alpha \circ \phi$ is $j^*$ above, so it is injective. This implies $\phi$ is injective. 

For the surjectivity, we already know the following:

    1. $K^*_G(G)$ is a free $R(G)$-module

    2. $K^*_T(G)=K^*_G(G) \otimes _R R(T) $
 
Now as $\alpha$ is injective,  Im$(\alpha )$ contains no torsion element.
Since 
    $\beta$ is the localization map, we see that the composition of $\beta \circ \alpha$ 
is injective.
 
Now by considering the Weyl group actions, we have:

Im$(\beta \circ \alpha) \subset K^*_T(T)^W$
and
$K^*_T(T)^W = {\rm Im}j^*$,
so we get the surjectivity. Q.E.D.

So we finished the proof of our main theorem.

\section{The proof of lemma 5.2}

First we prove a general result about holomorphic
differential forms $\Omega^i_{hol}(Y)$ on a complex manifold $Y$.

\begin {prop}
Let $X$ and $Y$ be  complex-analytic manifolds 
 and let $f:X\to Y$ be a holomorphic mapping.
Assume there is a finite group which acts on $X$
holomorphically, in such a way that $Y$ identifies with the quotient
of $X$ by $G$. Then we have
$$\Omega^i_{hol}(Y)=\Omega^i_{hol}(X)^G.$$
\end{prop}

\proof Let $\omega$ be a differential form on $X$ which
is $G$-invariant. Then $\omega$ is a meromorphic differential
form on $Y$, and we have to show it has no pole along any irreducible
divisor $D$ of $Y$. Pick a component $\Delta$ of $f^{-1}(D)$.
Near a general point $p$ of $\Delta$, we can find holomorphic
coordinates $(z_1,\cdots,z_n)$ on $X$ so that $\Delta$ has equation
$z_1=0$ and that $(z_1^e,\cdots,z_n)$ are holomorphic
coordinates on $Y$ near $f(p)\in Y$. We expand $\omega$ as

$$\omega=\sum_{I\subset \{ 2,\cdots,n\}}f_Idz_I
+dz_1\wedge\sum_{J\subset \{ 2,\cdots,n\}}g_Jdz_J.$$
Then the condition that $\omega$ is invariant under the transformation
$T(z_1,\cdots,z_n)=(\zeta z_1,z_2,\cdots,z_n)$ where $\zeta=e^{2\pi i\over e}$ translates
into the conditions

(1) $f_I$ is $T$-invariant

(2) $g_I(Tz)=\zeta^{-1}g_I(z)$.

The first condition implies that each term $f_Idz_I$ is a holomorphic
differential form in a neighborhood of $p$. The second condition
implies that $g_J$ is divisible by $z_1^{e-1}$, and therefore
is of the form $g_J=z_1^{e-1}h_J$, where $h_J$ is $T$-invariant, hence
descends to a holomorphic function in a neighborhood of $p$ in $Y$.
Then we have $g_Jdz_1\wedge dz_J={1\over e}h_Jd(z_1^e) \wedge dz_J$.
This concludes the proof. Q.E.D.

\begin {prop} Let ${\tt}$ be a split torus over $Spec(\BZ)$.
Let $X({\tt})$ denote the character group of ${\tt}$, and let $W$
be a finite subgroup of $Aut(X({\tt}))$. We can then view as a finite group
acting on ${\tt}$. Assume that the quotient scheme ${\tt}/W$ is smooth over
$Spec(\BZ)$. Then we have
$$\Omega^i(\tt)^W=\Omega^i(\tt/W).$$
\end{prop}
\proof Consider the projection map $\pi:\tt\to \tt/W$.
Let $\un\Omega^i_\tt$ denote the sheaf of Grothendieck differential
forms on $\tt$. Then we will show that the $W$-invariant part
$[\pi_*\un\Omega^i_\tt]^W$ of the direct image sheaf $\pi_*\un\Omega^i_\tt$
identifies with $\un\Omega^i_{\tt/W}$. There is an obvious
inclusion $$u:\un\Omega^i_{\tt/W}\hook [\pi_*\un\Omega^i_\tt]^W$$
As $\tt/W$ is smooth, the sheaf $\un\Omega^i_{\tt/W}$ is locally free and $\tt/W$
is normal. Hence it is enough to show that $u$ is an isomorphism
outside of a subscheme of codimension $\geq 2$.

We first construct a $W$-invariant closed subscheme $D$ of $\tt$
such that the restriction of $\pi$ to $U=\tt \setminus D$ is an etale
map $U\to U/W$. For $g\in W$, denote by $\tt ^g\subset \tt$ the fixed point scheme
of $g$, which is a closed group subscheme. Then we set
$$D=\cup_{g\in W,g\neq 1}~
\tt ^g.$$
Let $\Delta$ be the image of $D$ in $\tt /W$. As $\pi$ is proper, $\Delta$
is a closed subscheme and is the union of the $\Delta_g=\pi(\tt ^g)$.
As the map $U\to U/W$ is finite etale and Galois with Galois group $W$, it follows that the restriction
of $u$ to $U/W$ is an isomorphism. It is then enough to verify
that $u$ induces an isomorphism near a generic point of a component
$C$ of $\Delta$. Then $C$ is a component of $\Delta_g$ for some $g\in W$.
We only need to consider $g$ such that $\tt ^g$ is of codimension $1$. This
is equivalent to the condition
that $w$, as an automorphism of $X(\tt)$, is such that
$w-Id$ has rank $n-1$, where $n=dim(\tt)$. So if $w\neq 1$, the
only eigenvalues of $w$ are $1$ with multiplicity $n-1$ and $-1$
with multiplicity $1$. Then there are 2 possibilities:

(a) in a suitable basis $(\chi_1,\cdots,\chi_n)$ of $X(\tt)$,
we have $w\chi_1=-\chi_1$ and $w\chi_j=\chi_j$ for $j\geq 2$.

(b) in a suitable basis $(\chi_1,\cdots,\chi_n)$ of $X(\tt)$,
we have $w\chi_1=\chi_2$, $w\chi_2=\chi_1$ and $w\chi_j=\chi_j$
for $j\geq 3$.

In case (a), $\tt ^g$ is  the connected subscheme of $\tt$ of equation
$\chi^2=1$.  $\tt ^g$ has 2 irreducible
components, of equations $\chi_1=\pm 1$. These two components meet
in characteristic $2$. Note that each component of $\tt ^g$ maps
onto $Spec(\BZ)$. In case (b), $\tt ^g$ is irreducible and maps
onto $Spec(\BZ)$. Hence the same is true for $\Delta_g=\pi(\tt ^g)$.

Since every component of $\Delta$ which has codimension $1$ maps
onto $Spec(\BZ)$, it follows that its generic point $\eta$
maps to the generic point $Spec(\BQ)$ of $Spec(\BZ)$.
Hence to verify that $u$ is an isomorphism near $\eta$,
there is no harm in replacing $\tt$ with the algebraic torus
$\tt _{\BQ}=\tt \times_{Spec(\BZ)}Spec(\BQ)$ over $\BQ$. Then as $\BC$ is a field
extension of $\BQ$, we may work instead with 
$\tt _{\BC}=\tt \times_{Spec(\BZ)}Spec(\BC)$. Then we are reduced to showing
that if $V$ is a $W$-invariant open set (for the transcendental
topology) in $\tt _{\BC}$, then any
$W$-invariant Grothendieck differential $\omega$ on $V$ descends
to a Grothendieck differential on $V/W$. It is enough to show that
$\omega$ defines a holomorphic differential form on $V/W$, which follows
from Proposition 6.1. Q.E.D.

\noindent {\it Proof of lemma 5.2} We note that Proposition 6.2 implies lemma 5.2 easily.

For $G$, $T$, we have their complexification $G_{\C}$ and $T_ {\C}$, and $R(G_{\C})=R(G)$, $R(T_{\C})=R(T)$. For the torus $\tt =Spec (R(T_{\C}))$ over $Spec (\Z)$, ${\cal O}({\cal T})=R(T)$. 

We have $$\Omega ^*_{R(T)/\Z}= \Omega ^*_{{\cal T}/Spec(\Z)}$$
hence
$$ \begin {array}{ll}[\Omega ^*_{R(T)/\Z}]^W 
&=\Omega ^*_{({\cal T}/W)/Spec(\Z)}\\
&=\Omega ^*_{{\cal O}({\cal T}/W)/Spec(\Z)}\\
&=\Omega ^*_{R(T)^W/Spec(\Z)}\\
&=\Omega ^*_{R(G)/\Z}
\end{array}
$$
Q.E.D.

\section { The Case of $G=PSU(3)$}

So far we discussed the case when $G$ is simply connected. What happens when $G$ is non-simply connected? Let us consider $PSU(3)$.

Since this is a long computation, we just list the results here. $PSU(3)$ has a maximal torus $T^2$. There are elements $X_1$ and $X_2$  in $R(T^2)$:

$$\begin{array}{l}
X_1({\rm diag} (\l _1, \l _2, \l _3))=\frac {\l _2}{\l _1}\\
\\
X_2({\rm diag} (\l _1, \l _2, \l _3))=\l _1^3
\end{array}$$

We have 
$$R(T^2)=\Z[X_1, X_2, X_1^{-1}, X_2^{-1}]$$

By taking the Weyl group invariants of $R(T^2)$, we can get

\begin{lemma} $R(PSU(3))=\Z[Z_1, Z_2, Z_3]
\subset R(T^2)=\Z[X_1, X_2, X_1^{-1}, X_2^{-1}]$, where
$$\begin{array}{l}
Z_1=X_1+X_1^{-1}+X_1X_2+X_1^{-1} X_2^{-1}+X_1^2X_2+X_1^{-2}X_2^{-1}\\
Z_2=X_2+X_1^3X_2+X_1^{-3}X_2^{-2}\\
Z_3=X_2^{-1}+X_1^{-3}X_2^{-1}+X_1^3X_2^2
\end{array}$$
\end {lemma}

We see in fact that $R(PSU(3))$ is isomorphic to 
$\Z[X, Y, Z]/(X^3-YZ)$, so using the result of \cite {ll} or \cite {lll}, we obtain 

\begin{prop}: The Hochschild homology of $R=R(PSU(3))$ is given by
$$\begin{array}{l}
HH_0(R)=R\\
HH_1(R)=\Omega^1_{R/\Z}\\
HH_2(R)=\Omega^2_{R/\Z}\\
HH_3(R)=\Omega^3_{R/\Z}\cong T(\Omega^2_{R/\Z})\\
HH_k(R)=\Omega^3_{R/\Z}, \ k\ge 4
\end{array}
$$
\end{prop}

This is very interesting: the fact that $R$ has infinite homological dimension implies  the non-smoothness of 
$R(PSU(3))$. It is easy to see that the algebraic variety Spec($R\otimes \C$) is not smooth, as it 
has a singular point, corresponding to $A$= diag$(1, \omg, \omg ^2)\in PSU(3)$.

$K^0_{PSU(3)}(PSU(3))$ has a very nice $R(PSU(3))$-torsion part, it comes from following line bundle:
\begin{lemma} There exists a $PSU(3)$-equivariant line bundle $L$ over $PSU(3)$, such that:

$$\begin{array}{l}
1) L^{\otimes 3}=1\\
2) L  {\rm\ is\ trivial\ outside\ the\ orbit\ G.A}\\
3) R(PSU(3))/Ann([L]-1)=\Z \\
4) p_*1=1\oplus L\oplus L^2
\end{array}
$$
Here $[L]$ denotes the class of $L$ in $K_{PSU(3)}^*(PSU(3))$, $p:\ SU(3)\to PSU(3)$ is the quotient map.
\end{lemma}

By using the affine Weyl group \cite {jh1}, we get the orbit space as in figure 1, where points 
in the edges with arrows are identified by $(a,b)\tl (a, 1)$, 
the orbit passing through $A$ corresponds to $=(-\frac 23, 1)$.
 
\begin{figure}[ht]
\psfig{figure=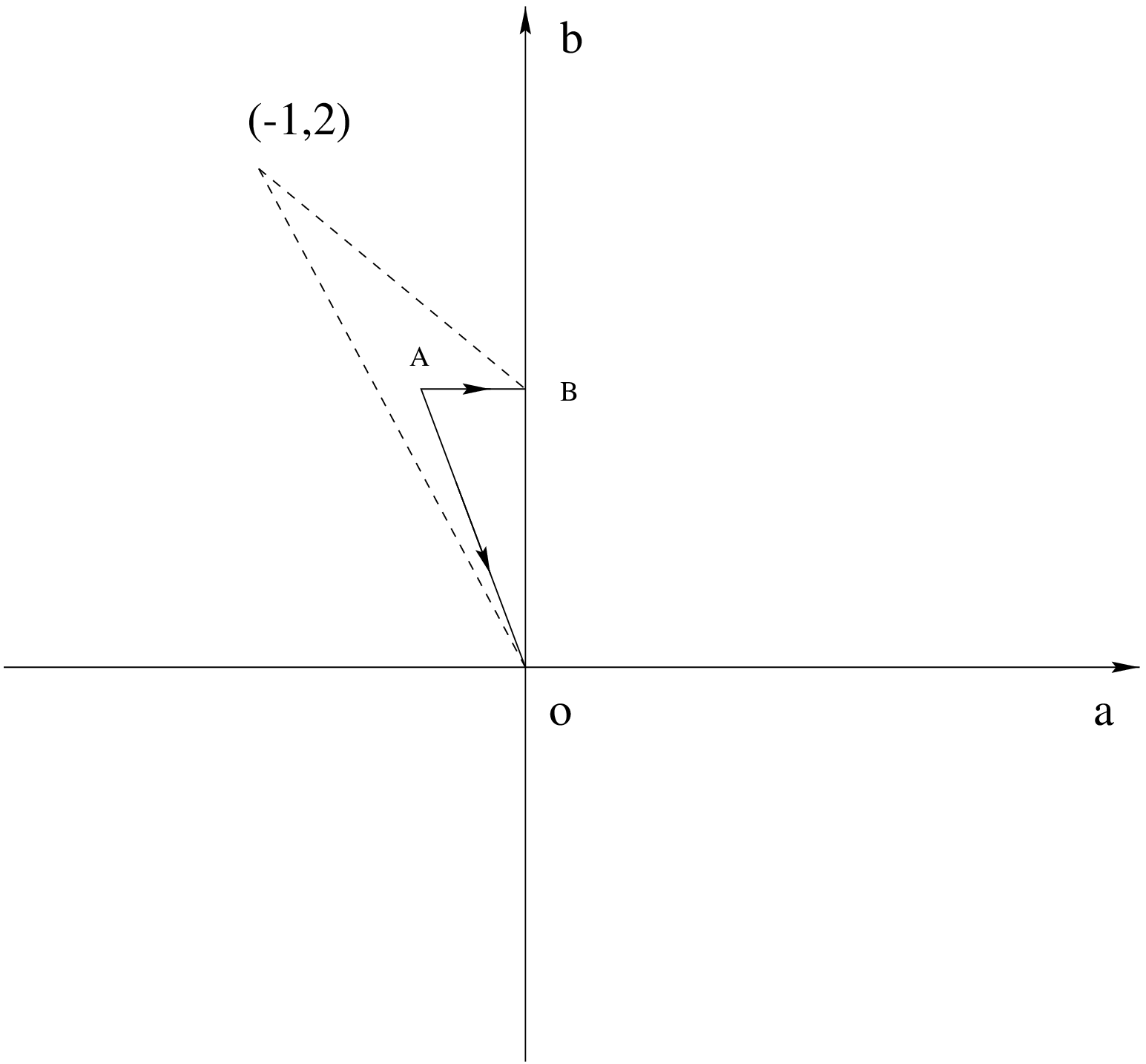}
\caption {Orbit Space of PSU(3)}
\end{figure}

Now we can construct a very nice covering of $PSU(3)$, then by the Segal spectral sequence, we find that:

\begin{prop}
$$K_{PSU(3)}^0(PSU(3))=R(PSU(3))^2 \oplus \Z .([L]-1)\oplus \Z .([L]^2-1)$$
While $K_{PSU(3)}^1(PSU(3))$ is torsion free.
\end{prop}

So in this case, we have $\Omega ^*_{R/\Z}\not \cong K^*_G(G)$. The general case  remains an open problem.

\end{document}